\documentstyle[amscd,amssymb,verbatim,12pt]{amsart}
\newcommand{\Aff}{\mbox{\rm Aff}}

\newcommand{\Aut}{\mbox{\rm Aut}}

\newcommand{\C}{{\Bbb C}}

\newcommand{\const}{\mbox{\rm const.}}

\newcommand{\diam}{\mbox{\rm diam}}
\newcommand{\Diff}{\mbox{\rm Diff}}

\newcommand{\End}{\mbox{\rm End}}

\newcommand{\HH}{\mbox{\rm H}}

\newcommand{\Id}{\mbox{\rm Id}}

\newcommand{\Image}{\mbox{\rm Im}}

\newcommand{\inj}{\mbox{\rm inj}}

\newcommand{\R}{{\Bbb R}}

\newcommand{\SL}{\mbox{\rm SL}}
\newcommand{\SO}{\mbox{\rm SO}}
\newcommand{\spann}{\mbox{\rm span}}

\newcommand{\Spin}{\mbox{\rm Spin}}

\newcommand{\Tr}{\mbox{\rm Tr}}

\newcommand{\vol}{\mbox{\rm vol}}
\newcommand{\Z}{{\Bbb Z}}
\theoremstyle{plain}
\newtheorem{definition}{Definition}

\newtheorem{lemma}{Lemma}
\newtheorem{theorem}{Theorem}
\newtheorem{proposition}{Proposition}

\numberwithin{equation}{section}
\renewcommand{\rm}{\normalshape}
\begin{document}
\title{Collapsing and Dirac-Type Operators}
\author{John Lott}
\address{Department of Mathematics\\
University of Michigan\\
Ann Arbor, MI  48109-1109\\
USA}
\email{lott@@math.lsa.umich.edu}
\thanks{Research supported by NSF grant DMS-9704633}
\date{May 1, 2000}
\maketitle
\begin{abstract}
We analyze the limit of the spectrum of a geometric Dirac-type
operator under a collapse with bounded
diameter and bounded sectional curvature. 
In the case of a smooth limit space $B$, 
we show that the limit of the spectrum is given by the spectrum of
a certain first-order differential operator on $B$, which can be 
constructed using superconnections.
In the case of a general limit space $X$, we express the limit operator in
terms of
a transversally elliptic operator on a $G$-manifold $\check{X}$ with
$X \: = \: \check{X}/G$.
As an application, we give a characterization of manifolds which do not admit
uniform upper
bounds, in terms of diameter and
sectional curvature, on the $k$-th eigenvalue of the square
of a Dirac-type operator. We also give a formula for the
essential spectrum of 
a Dirac-type operator on a finite-volume manifold
with pinched negative sectional curvature.
\end{abstract}

\section{Introduction} \label{sect1}

In a previous paper we analyzed the limit of the spectrum of 
the differential form Laplacian on a manifold,
under a collapse with bounded
diameter and bounded sectional curvature \cite{Lott (1999)}. 
In the present paper, we extend the analysis
of \cite{Lott (1999)}
to geometric Dirac-type operators. As the present paper is a sequel to
\cite{Lott (1999)}, we refer to the introduction of \cite{Lott (1999)} for
background information about collapsing with bounded curvature 
and its relation to analytic questions.

Let $M$ be a connected closed oriented Riemannian manifold of dimension
$n \: > \: 0$.
If $M$ is spin then we put $G \: = \: \Spin(n)$ and if $M$ is not spin then
we put $G \: = \: \SO(n)$. The spinor-type fields that we consider are 
sections of a vector bundle $E^M$
associated to a $G$-Clifford module $V$, the latter being 
in the sense of Definition \ref{def2} of Section \ref{sect2}. 
The ensuing Dirac-type operator $D^M$ acts on sections of $E^M$.
We will 
think of the spectrum $\sigma(D^M)$ of $D^M$ as a set of real numbers with
multiplicities, corresponding to possible multiple eigenvalues.
For simplicity, in this introduction we will sometimes refer to the
Dirac-type operators as acting on spinors, even though the results are
more general.

We first consider a collapse in which the limit space is a
smooth Riemannian manifold. The model case
is that of a Riemannian affine fiber bundle.

\begin{definition} \label{def1}
A Riemannian affine fiber bundle is a smooth fiber bundle 
$\pi \: : \: M \rightarrow B$ whose
fiber $Z$ is an infranilmanifold and whose structure group is reduced from
$\Diff(Z)$ to $\Aff(Z)$, along with
\begin{itemize}
\item A horizontal distribution $T^HM$ whose holonomy lies in $\Aff(Z)$,
\item A family of vertical Riemannian metrics $g^{TZ}$ which are parallel
with respect to the flat affine connections on the fibers $Z_b$ and
\item A Riemannian metric $g^{TB}$ on $B$.
\end{itemize}
\end{definition}

Given a Riemannian affine fiber bundle $\pi \: : \: M \rightarrow B$,
there is a Riemannian metric $g^{TM}$ on $M$ constructed from
$T^HM$, $g^{TZ}$ and $g^{TB}$.
Let $R^M$ denote the Riemann curvature tensor of $(M, g^{TM})$,
let $\Pi$ denote the second fundamental forms of the fibers 
$\{Z_b\}_{b \in B}$ and
let $T \in \Omega^2(M; TZ)$ be the curvature of $T^HM$.
Given $b \in B$, there is a natural flat connection on $E^M \big|_{Z_b}$
which is constructed using the affine structure of $Z_b$. We define a
Clifford bundle $E^B$ on $B$ whose fiber over $b \in B$ consists of
the parallel sections of $E^M \big|_{Z_b}$. The operator $D^M$ restricts
to a first-order differential operator $D^B$ on $C^\infty(B; E^B)$. 
If $V$ happens to be the spinor module then
we show that $D^B$ is the ``quantization'' of a certain superconnection
on $B$. For general $V$, there is an additional zeroth-order term in $D^B$
which depends on $\Pi$ and $T$.

We show that the spectrum
of $D^M$ coincides with that of $D^B$ up
to a high level, which depends on the maximum diameter 
$\diam(Z)$ of the fibers $\{Z_b\}_{b \in B}$.

\begin{theorem} \label{thm1}
There are positive constants $A$, $A^\prime$ and $C$ which only depend on
$n$ and $V$ such that if $\parallel R^Z \parallel_\infty \:
\diam(Z)^2 \: \le \: A^\prime$ then the intersection of $\sigma(D^M)$ with
the interval
\begin{align} \label{1.1}
[ 
- \left( \: A \: \diam(Z)^{-2} \: - \: \right.
& 
\left. C \: \left( 
\parallel R^M \parallel_\infty \: + \: {\parallel \Pi \parallel_\infty}^2
\: + \: {\parallel T \parallel_\infty}^2 \right) \right)^{1/2},  \notag  \\
&  
\left( A \: \diam(Z)^{-2} \: - \: C \: \left( 
\parallel R^M \parallel_\infty \: + \: {\parallel \Pi \parallel_\infty}^2
\: + \: {\parallel T \parallel_\infty}^2 \right) \right)^{1/2}
] 
\end{align}
equals the intersection of $\sigma(D^B)$ with (\ref{1.1}).
\end{theorem} 

If $Z \: = \: S^1$, $\Pi \: = \: 0$ and $V$ is the spinor module then we 
recover some results of
\cite[Section 4]{Ammann-Bar (1998)}; see also
\cite[Theorem 1.5]{Dai (1991)}. The proof of Theorem \ref{thm1} follows
the same strategy as the proof of the analogous 
\cite[Theorem 1]{Lott (1999)}. Consequently, in the proof of Theorem 
\ref{thm1},
we only indicate the changes that need to be made in the proof of 
\cite[Theorem 1]{Lott (1999)} and refer to \cite{Lott (1999)} for details.

Given $B$, Cheeger, Fukaya and Gromov showed that under some
curvature bounds, any Riemannian
manifold $M$ which is sufficiently
Gromov-Hausdorff close to $B$ can be well approximated by a Riemannian
affine fiber bundle \cite{Cheeger-Fukaya-Gromov (1992)}.
Using this fact, we show that the spectrum of $D^M$ can be 
uniformly approximated 
by that of a certain first-order differential
operator $D^B$ on $B$, at least up to a high level 
which depends
on the Gromov-Hausdorff distance between $M$ and $B$.

Given $\epsilon > 0$ and two collections of
real numbers $\{a_i\}_{i \in I}$ and $\{b_j\}_{j \in J}$,
we say that
$\{a_i\}_{i \in I}$ and $\{b_j\}_{j \in J}$ are $\epsilon$-close if
there is a bijection $\alpha : I \rightarrow J$ such that for all
$i \in I$, $|b_{\alpha(i)} \: - \: a_i| \: \le \: \epsilon$. 

\begin{theorem} \label{thm2}
Let $B$ be a fixed smooth connected closed Riemannian manifold. Given
$n \in \Z^+$, take  $G \: \in \{ \SO(n), \Spin(n)\}$ and let $V$ be a 
$G$-Clifford module.
Then for any $\epsilon \: > \: 0$ and $K \: > \: 0$, there are positive
constants $A(B, n, V, \epsilon, K)$, $A^\prime(B, n, V, \epsilon, K)$, 
and $C(B, n, V, \epsilon, K)$
so that the following holds. 
Let $M$ be an
$n$-dimensional connected closed oriented Riemannian manifold with
a $G$-structure such that $\parallel R^M \parallel_\infty \: \le \: K$ and
$d_{GH}(M, B) \: \le \: A^\prime$. Then there are a Clifford module
$E^B$ on $B$ and a certain first-order differential operator
$D^B$ on $C^\infty(B; E^B)$ such that \\
1. $\left\{\sinh^{-1} \left( \frac{\lambda}{\sqrt{2K}} \right) 
\: : \: \lambda \in \sigma(D^M) \text{ and } 
\lambda^2 \: \le \:
A \: d_{GH}(M, B)^{-2} \: - \: C \right\}$ is $\epsilon$-close to a subset of
$\left\{\sinh^{-1} \left( \frac{\lambda}{\sqrt{2K}} \right) 
\: : \: \lambda \in \sigma(D^B) \right\}$, and \\
2. $\left\{\sinh^{-1} \left( \frac{\lambda}{\sqrt{2K}} \right) 
\: : \: \lambda \in \sigma(D^B) \text{ and } 
\lambda^2 \: \le \:
A \: d_{GH}(M, B)^{-2} \: - \: C \right\}$ is $\epsilon$-close to a subset of
$\left\{\sinh^{-1} \left( \frac{\lambda}{\sqrt{2K}} \right) 
\: : \: \lambda \in \sigma(D^M) \right\}$.
\end{theorem}

The other results in this paper concern collapsing to a possibly-singular
space.  Let $X$ be a limit space of a sequence $\{M_i\}_{i=1}^\infty$ 
of $n$-dimensional connected closed
oriented Riemannian manifolds with uniformly bounded diameter and
uniformly bounded sectional curvature.  In general, $X$ is not 
homeomorphic to
a manifold.  However, Fukaya showed that $X$ is homeomorphic to $\check{X}/G$,
where $\check{X}$ is a manifold and $G$ is a compact Lie group which acts
on $\check{X}$ \cite{Fukaya (1988)}.  This comes from
writing $M_i \: = \: P_i/G$, where $G \: = \: \SO(n)$ and
$P_i$ is the oriented orthonormal frame
bundle of $M_i$. There is a canonical Riemannian metric on $P_i$. Then
$\{P_i\}_{i=1}^\infty$ has a subsequence which Gromov-Hausdorff
converges to a manifold $\check{X}$. As the convergence argument can be
done $G$-equivariantly, the corresponding subsequence of $\{M_i\}_{i=1}^\infty$
converges to $X \: = \: \check{X}/G$. In general, $\check{X}$ is a
smooth manifold with a metric which is $C^{1,\alpha}$ regular for all
$\alpha \in (0,1)$.

In \cite{Lott (1999)} we dealt with the limit of the spectra of the
differential form
Laplacians $\{\triangle^{M_i}\}_{i=1}^\infty$ 
on the manifolds $\{M_i\}_{i=1}^\infty$. We defined a
limit operator $\triangle^X$ 
which acts on the ``differential forms'' on $X$, coupled to a 
superconnection. In order
to make this precise, we defined the ``differential forms'' on $X$ to be
the $G$-basic differential forms on $\check{X}$. We constructed the
corresponding differential form Laplacian $\triangle^X$ and showed 
that its spectrum described the limit of the spectra of 
$\{\triangle^{M_i}\}_{i=1}^\infty$. We refer to \cite{Lott (1999)} for
the precise statements.

In the case of geometric 
Dirac-type operators $D^{M_i}$, there is a fundamental 
problem in
extending this approach.
Namely, if $\check{X}$ is a spin manifold on
which a compact Lie group $G$ acts isometrically and preserving the spin 
structure then
there does not seem to be a notion of $G$-basic spinors on $\check{X}$. In
order to get around this problem, we take a different approach. For a given
$n$-dimensional Riemannian spin manifold $M$, 
put $G \: = \: \Spin(n)$, let $P$ be the principal $\Spin(n)$-bundle
of $M$ and let $V$ be the spinor module.
One can identify the spinor fields on $M$ with
$\left( C^\infty(P) \: \otimes \: V \right)^G$, the $G$-invariant subspace
of $C^\infty(P) \: \otimes \: V$.
There are canonical horizontal vector fields 
$\{{\frak Y}_j\}_{j=1}^n$ on $P$ and the Dirac operator takes the form
$D^{M} \: = \: - \: i \: \sum_{j=1}^n \gamma^j \: {\frak Y}_j$.
Furthermore, $(D^M)^2$ can be written in a particularly simple form.
As in equation (\ref{4.2}) below, when acting on
$\left( C^\infty(P) \otimes V \right)^G$, $(D^M)^2$ becomes the
scalar Laplacian on $P$ (acting on $V$-valued functions) plus a 
zeroth-order term.

Following this viewpoint, 
it makes sense to define the limiting ``spinor fields'' on $X$ to be the
elements of
$\left( C^\infty(\check{X}) \otimes V \right)^G$. 
We can then extend Theorem \ref{thm1} to the setting of $G$-equivariant
Riemannian affine fiber bundles. Namely, the limit operator
$D^X$ turns out to be a $G$-invariant first-order differential operator on 
$C^\infty(\check{X}) \otimes V$, transversally elliptic 
in the sense of Atiyah \cite{Atiyah (1974)}, which one then restricts
to the $G$-invariant subspace $\left( C^\infty(\check{X}) \otimes V \right)^G$.
In Theorem \ref{thm6} below, 
we show that the analog of Theorem \ref{thm1} holds, in which
$D^B$ is replaced by $D^X$.

Theorem \ref{thm6} refers to a given $G$-equivariant Riemannian affine
fiber bundle.
In order to deal with arbitrary collapsing sequences,
we use the aforementioned representation of $(D^M)^2$ as a Laplace-type
operator on $P$. If $\{M_i\}_{i=1}^\infty$ is a sequence of
$n$-dimensional
Riemannian manifolds with uniformly bounded diameter
and uniformly bounded sectional curvature then we show that after taking
a subsequence, the spectra of $\{(D^{M_i})^2\}_{i=1}^\infty$ 
converge to the spectrum of a
Laplace-type operator on a limit space. Let 
$\{\lambda_k(|D^M|)\}_{k=1}^\infty$ denote
the eigenvalues of $|D^M|$, counted with multiplicity.

\begin{theorem} \label{thm3}
Given $n \in \Z^+$ and $G \in \{\SO(n), \Spin(n)\}$,
let $\{M_i\}_{i=1}^\infty$ be a sequence of connected closed 
oriented $n$-dimensional
Riemannian manifolds with a $G$-structure.
Let $V$ be a $G$-Clifford module.
Suppose that for some $D, \: K \: > \: 0$ and
for each $i \in \Z^+$, we have $\diam(M_i) \: \le \: D$ and
$\parallel R^{M_i} \parallel_\infty \: \le \: K$.
Then there are\\
1. A subsequence of $\{M_i\}_{i=1}^\infty$, which we relabel as
$\{M_i\}_{i=1}^\infty$,\\
2. A smooth closed $G$-manifold $\check{X}$ with a $G$-invariant
Riemannian metric 
$g^{T\check{X}}$ which is
$C^{1,\alpha}$-regular for all $\alpha \in (0,1)$,\\
3. A positive $G$-invariant function $\chi \in C(\check{X})$ with
$\int_{\check{X}} \chi \: d\vol \: = \: 1$ and\\
4. A $G$-invariant function ${\cal V} \in L^\infty(\check{X}) 
\otimes \End(V)$\\
such that if $\triangle^{\check{X}}$ denotes the Laplacian on
$L^2(\check{X}, \chi \: d\vol) \otimes V$ \cite[(0.8)]{Fukaya (1987)} and 
$|D^X|$ denotes the operator
$\sqrt{\triangle^{\check{X}} \: + \: {\cal V}}$ acting on 
$\left( L^2(\check{X}, \chi \: d\vol) \otimes V \right)^G$ then for all 
$k \in \Z^+$,
\begin{equation} \label{1.4}
\lim_{i \rightarrow \infty} \lambda_k \left( |D^{M_i}| \right) \: = \:
\lambda_k \left( |D^{X}| \right).
\end{equation}
\end{theorem}

In the special case of the signature operator,
the proof of Theorem \ref{thm3} is somewhat simpler than that of the analogous
\cite[Proposition 11]{Lott (1999)}, in that we essentially
only have to deal with scalar Laplacians.  
However, \cite[Proposition 11]{Lott (1999)} gives
more detailed information. In particular, it expresses the limit
operator in terms
of a basic flat degree-$1$ superconnection on $\check{X}$.
This seems to be
necessary in order to prove the results of \cite{Lott (1999)} concerning small
eigenvalues. Of course, one does not expect to have analogous results 
concerning the small eigenvalues of general geometric Dirac-type operators,
as their zero-eigenvalues have no topological
meaning.

As an application of Theorem \ref{thm3}, we give a characterization of
manifolds which do not have a uniform upper bound on the $k$-th 
eigenvalue of $|D^M|$, in terms of diameter and sectional curvature.

\begin{theorem} \label{thm4}
Let $M$ be a connected closed oriented manifold with a $G$-structure. 
Let $V$ be a
$G$-Clifford module.  Suppose that 
for some $K \: > \: 0$ and $k \in \Z^+$, 
there is no uniform upper bound on $\lambda_k(|D^M|)$ among Riemannian
metrics on $M$ with $\diam(M) \: = \: 1$ and 
$\parallel R^M \parallel_\infty \: \le \: K$. Then $M$ admits a 
possibly-singular fibration $M \rightarrow X$ whose generic fiber is an
infranilmanifold $Z$ such that the restriction of
$E^M$ to $Z$ does not have any nonzero 
affine-parallel sections.
\end{theorem}

More precisely, the possibly-singular fibration $M \rightarrow X$
of Theorem \ref{thm4} is
the $G$-quotient of a $G$-equivariant Riemannian affine fiber bundle
$P \rightarrow \check{X}$.
Theorem \ref{thm4} is an analog of \cite[Theorem 5.2]{Lott (1999)}.
A simple example of Theorem \ref{thm4} comes from considering spinors on
$M \: = \: S^1 \times N$, where $N$ is a spin manifold and the spin
structure on $S^1$ is the one that does not admit a harmonic spinor.
Upon shrinking the $S^1$-fiber, the eigenvalues of $D_M$ go off to 
$\pm \: \infty$. More generally, let $\pi : M \rightarrow B$ be
an affine fiber bundle.
Theorem \ref{thm1} implies that if 
$E^M \big|_Z$ does not have any nonzero
affine-parallel sections then upon 
collapsing $M$ to $B$ as in
\cite[Section 6]{Fukaya (1989)}, the eigenvalues of $D_M$ go off to
$\pm \: \infty$.

Finally, we give a result about the essential spectrum of a geometric
Dirac-type operator on a finite-volume manifold of pinched negative 
curvature, which is an analog of
\cite[Theorem 2]{Lott2 (1999)}.
Let $M$ be a complete connected oriented $n$-dimensional 
Riemannian manifold with a 
$G$-structure. Suppose that $M$ has finite volume and its 
sectional curvatures satisfy $- \: b^2 \: \le \: K
\: \le \: - \: a^2$, with $0 \: < \: a \: \le \: b$. Let $V$ be a $G$-Clifford
module. Label the ends of
$M$ by $I \in \{1,\ldots, N\}$. An end of $M$ has a neighborhood $U_I$
whose closure is homeomorphic to $[0, \infty) \: \times \: Z_I$, where 
the first coordinate is the Busemann function corresponding to a ray
exiting the end, and $Z_I$ is an infranilmanifold. 
Let $E^M$ be the vector bundle on $M$ associated to the pair $(G, V)$
and let $D^M$ be the corresponding Dirac-type operator.
If $U_I$ lies far enough out the end then
for each $s \in [0, \infty)$,
$C^\infty \left( \{s\} \times Z_I; E^M \big|_{\{s\} \times Z_I} \right)$ 
decomposes as the direct sum of a 
finite-dimensional space $E^B_{I,s}$, consisting of ``bounded energy''
sections, and its orthogonal complement, consisting of ``high energy''
sections.
The vector spaces $\{ E^B_{I,s} \}_{s \in [0, \infty)}$ fit
together to form a vector bundle $E^B_I$ on $[0, \infty)$. Let $P_0$
be orthogonal projection from $\bigoplus_{I=1}^N C^\infty \left(\overline{U_I};
E^M \big|_{\overline{U_I}} \right)$ to 
$\bigoplus_{I=1}^N C^\infty \left([0, \infty); E^B_I \right)$. 
Let $D^M_{end}$ be the
restriction of $D^M$ to $\bigoplus_{I=1}^N C^\infty \left(\overline{U_I};
E^M \big|_{\overline{U_I}}\right)$, 
say with Atiyah-Patodi-Singer boundary conditions.
Then $P_0 \:  D^M_{end} \: P_0$ is a first-order ordinary differential
operator on $\bigoplus_{I=1}^N C^\infty \left( [0, \infty); E^B_I \right)$.

\begin{theorem} \label{thm5}
The essential spectrum of $D^M$ is the same as that of 
$P_0 \:  D^M_{end} \: P_0$.
\end{theorem}

There is some intersection between Theorem \ref{thm5} and the results of
\cite[Theorem 0.1]{Ballmann-Bruning (1999)}, concerning the essential
spectrum of $D^M$ when $n \: = \: 2$ and under an additional curvature
assumption, and \cite[Theorem 1]{Bar (1998)}, concerning the
essential spectrum of $D^M$ when $M$ is hyperbolic and $V$ is the
spinor module.

I thank the Max-Planck-Institut-Bonn for its hospitality while this
research was performed.

\section{Dirac-type Operators and Infranilmanifolds} \label{sect2}
Given $n \in \Z^+$, let $G$ be either $\SO(n)$ or $\Spin(n)$.

\begin{definition} \label{def2}
A $G$-Clifford module consists of 
a finite-dimensional Hermitian $G$-vector space $V$ and
a $G$-equivariant linear map $\gamma \: : \: \R^n \rightarrow \End(V)$ such
that $\gamma(v)^2 \: = \: |v|^2 \: \Id.$ 
and $\gamma(v)^* \: = \: \gamma(v)$. 
\end{definition}

Let $M$ be a connected closed oriented smooth $n$-dimensional Riemannian
manifold.
Put $G \: = \: \Spin(n)$ or $G \: = \: \SO(n)$, according as to whether or
not $M$ is spin. If $M$ is spin, fix a spin structure. Let $P$ be the 
corresponding
principal $G$-bundle, covering the oriented orthonormal frame bundle. 
Its topological isomorphism class is independent of the choice of 
Riemannian metric. Given the Riemannian metric,
there is a canonical $\R^n$-valued $1$-form $\theta$ on $P$, the
soldering form.

With respect to the standard basis $\{e_j\}_{j=1}^n$ of $\R^n$, we write
$\gamma^j \: = \: \gamma(e_j)$. We also take a basis
$\{\sigma^{ab}\}_{a,b=1}^n$ for the representation 
of the Lie algebra ${\frak g}$ on
$V$, so that $\left( \sigma^{ab} \right)^* \: = \: - \: \sigma^{ab}$ and
\begin{equation} \label{2.1}
[\sigma^{ab}, \sigma^{cd}] \: = \: 
\delta^{ad} \: \sigma^{bc} \: - \: \delta^{ac} \: \sigma^{bd} \: + \: 
\delta^{bc} \: \sigma^{ad} \: - \: \delta^{bd} \: \sigma^{ac}.
\end{equation}
The $G$-equivariance of $\gamma$ implies
\begin{equation} \label{2.2}
[ \gamma^a, \: \sigma^{bc}] \: = \:
\delta^{ab} \: \gamma^c \: - \: \delta^{ac} \: \gamma^b.
\end{equation}
{\bf Examples : } \\
1. If $G \: = \: \Spin(n)$ and $V$ is the spinor representation of $G$ then
$\sigma^{ab} \: = \: \frac{1}{4} \: [\gamma^a, \gamma^b]$. \\
2. If $G \: = \: \SO(n)$ and $V \: = \: \Lambda^*(\R^n) \otimes_\R \C$, let
$E^j$ and $I^j$ denote exterior and interior multiplication by $e^j$,
respectively. Put $\gamma^j \: = \: i \: \left( E^j \: - \: I^j \right)$ and
$\widehat{\gamma}^j \: = \: E^j \: + \: I^j$. Then
$\sigma^{ab} \: = \: \frac{1}{4} \: \left( [\gamma^a, \gamma^b] \: + \:
[\widehat{\gamma}^a, \widehat{\gamma}^b] \right)$. \\

Put $E^M \: = \: P \times_G V$. 
The Dirac-type operator $D^M$ acts on the space $C^\infty(M; E^M)$. As the
topological vector space
$C^\infty(M; E^M)$ is independent of any choice of Riemannian metric on
$M$, it makes sense to compare Dirac-type operators for different Riemannian
metrics on $M$; see \cite[Section 2]{Lott (1997)} for further discussion.

Let $g^{TM}$ be the Riemannian metric on $M$. Let $\omega$ be the Levi-Civita
connection on $P$. 
Let $\{e_j\}_{j=1}^n$ be a local oriented orthonormal basis
of $TM$, with dual basis $\{\tau^j\}_{j=1}^n$. Then 
we can write $\omega$ locally as a matrix-valued $1$-form
$\omega^a_{\: \: b} \: = \: \sum_{j=1}^n \omega^a_{\: \: bj} \tau^j$, and
\begin{equation} \label{2.3}
D^M \: = \: - \: i \: \sum_{j=1}^n \gamma^j \: \nabla_{e_j} \: = \:
- \: i \: \sum_{j=1}^n \gamma^j \: \left( e_j \: + \: \frac{1}{2} \: 
\sum_{a,b=1}^n
\omega_{abj} \: \sigma^{ab} \right).
\end{equation}
We have the Bochner-type equation
\begin{equation} \label{2.4}
(D^M)^2 \: = \: \nabla^* \nabla \: - \: \frac{1}{8} \: \sum_{a,b,i,j=1}^n
R^M_{abij} \: (\gamma^i \: \gamma^j \: - \: \gamma^j \: \gamma^i) \:
\sigma^{ab}. 
\end{equation}

As the set of Riemannian metrics on $M$ is an affine
space modeled on a Fr\'echet space, it makes sense to talk about an
analytic $1$-parameter family $\{c(t)\}_{t \in [0,1]}$ of metrics.
Then for $t \in [0,1]$, $\dot{c}(t)$ is a symmetric
$2$-tensor on $M$. Let $\parallel \dot{c}(t) 
\parallel_{c(t)}$ denote
the norm of $\dot{c}(t)$ with respect to $c(t)$, i.e.
\begin{equation} \label{2.5}
\parallel \dot{c}(t) 
\parallel_{c(t)} \: = \: \sup_{v \in TM - 0}
\frac{|\dot{c}(t)(v,v)|}{c(t)(v,v)}.
\end{equation}
Put $l(c) \: = \: \int_0^1 \parallel \dot{c}(t) 
\parallel_{c(t)} \: dt$. We extend the definition of $l(c)$ to 
piecewise-analytic families of metrics in the obvious way.
Given $K \: > \: 0$, 
let ${\cal M}(M, K)$ be the set of Riemannian metrics on $M$ with
$\parallel R^M \parallel_\infty 
\: \le \: K$.
Let $d$ be the corresponding length metric on ${\cal M}(M, K)$, computed
using piecewise-analytic paths in ${\cal M}(M, K)$.
Let $\sigma(D^M, g^{TM})$ 
denote the spectrum of $D^M$ as computed with $g^{TM}$, a discrete subset
of $\R$ which is counted with
multiplicity.

\begin{proposition} \label{prop1}
There is a constant $C \: = C(n, V) \: > \: 0$ such that for all
$K \: > \: 0$ and $g_1^{TM}, g_2^{TM} \: \in \: {\cal M}(M, K)$,
\begin{equation} \label{2.6}
\left\{ \sinh^{-1} \left( \frac{\lambda}{\sqrt{K}} \right) \: : \:
\lambda \in \sigma(D^M, g_1^{TM}) \right\}
\end{equation}
 and
\begin{equation} \label{2.7}
\left\{ \sinh^{-1} \left( \frac{\lambda}{\sqrt{K}} \right) \: : \:
\lambda \in \sigma(D^M, g_2^{TM}) \right\}
\end{equation}
are $C \: d(g_1^{TM}, g_2^{TM})$-close.
\end{proposition}
\begin{pf}
It is enough to show that there is a number $C$ such that 
if $\{c(t)\}_{t \in [0,1]}$ is an analytic $1$-parameter family of metrics
contained in ${\cal M}(M, K)$ then 
$\left\{ \sinh^{-1} \left( \frac{\lambda}{\sqrt{K}} \right) \: : \:
\lambda \in \sigma(D^M, c(0)) \right\}$ and
$\left\{ \sinh^{-1} \left( \frac{\lambda}{\sqrt{K}} \right) \: : \:
\lambda \in \sigma(D^M, c(1)) \right\}$ are
$C \: d(c(0), c(1))$-close.
By eigenvalue perturbation theory \cite[Chapter XII]{Reed-Simon (1978)}, the
subset $\bigcup_{t \in [0,1]} \{t\} \times \sigma(D^M, c(t))$ of $\R^2$ is the
union of the graphs of functions $\{\lambda_j(t)\}_{j \in \Z}$
which are analytic in $t$. Thus it is enough to show that for each $j \in \Z$,
\begin{equation} \label{2.8}
\left| \sinh^{-1} \left( \frac{\lambda_j(1)}{\sqrt{K}} \right) \: - 
\sinh^{-1} \left( \frac{\lambda_j(0)}{\sqrt{K}} \right) \right| \: \le \:
C \: l(c).
\end{equation}

Let $D(t)$ denote the Dirac-type operator constructed with the metric
$c(t)$. It is self-adjoint when acting on $L^2(E^M, d\vol(t))$. In order
to have all of the operators $\{D(t)\}_{t \in [0,1]}$ acting on the same
Hilbert space, define $f(t) \in C^\infty(M)$ by $f(t) \: = \:
\frac{d vol(t)}{d vol(0)}$. Then the spectrum of $D(t)$, acting
on $L^2(E^M, d\vol(t))$, is the same as the spectrum of the
self-adjoint operator $f(t)^{1/2} \: D(t) \: f(t)^{- \: 1/2}$ 
acting on
$L^2(E^M, d\vol(0))$.
One can now compute $\frac{d\lambda_j}{dt}$ using
eigenvalue perturbation theory, as in \cite[Chapter XII]{Reed-Simon (1978)}.
Let $\psi_j(t)$ be a smoothly-varying unit eigenvector whose eigenvalue is
$\lambda_j(t)$.
Define a quadratic form $T(t)$ on $TM$ by
\begin{equation} \label{2.9}
T(t)(X,Y) \: = \: \langle \psi_j, - \: i \: (\gamma(X) \nabla_Y \psi_j \: + \:
\gamma(Y) \nabla_X \psi_j ) \rangle \: + \: \langle
- \: i \: (\gamma(X) \nabla_Y \psi_j \: + \:
\gamma(Y) \nabla_X \psi_j ), \psi_j \rangle.  
\end{equation}
Using the metric 
$c(t)$ to convert the symmetric tensors $\dot{c}(t)$ and $T(t)$ to 
self-adjoint sections of  $\End(TM)$, one
finds 
\begin{equation} \label{2.10}
\frac{d\lambda_j}{dt} \: = \: - \: \frac{1}{8} \int_M \Tr \left(
\dot{c}(t) \: T(t) \right) \: d\vol(t).
\end{equation}
(This equation was shown for the pure Dirac operator, by different means,
in \cite{Bourguignon-Gauduchon (1992)}.) Then
\begin{equation} \label{2.11}
\left| \frac{d\lambda_j}{dt} \right| \: \le \: \const \:
\parallel \dot{c}(t) 
\parallel_{c(t)} \: \int_M \Tr(|T(t)|) \: d\vol(t).
\end{equation}
Letting $\{x_i\}_{i=1}^n$ be an orthonormal basis of eigenvectors of
$T(t)$ at a point $m \in M$, 
we have $\Tr(|T(t)|) \: = \: \sum_{i=1}^n |T(t)(x_i, x_i)|$. Then 
from (\ref{2.9}), we obtain
\begin{equation} \label{2.12}
\int_M \Tr(|T(t)|) \: d\vol(t) \: \le \: \const \left( \int_M |\nabla \psi_j|^2
\: d\vol(t) \right)^{1/2}.
\end{equation}
From (\ref{2.4}),
\begin{equation} \label{2.13}
\int_M |\nabla \psi_j|^2
\: d\vol(t) \: \le \: \lambda_j^2 \: + \: \const \: K.
\end{equation}
In summary, from (\ref{2.11}), (\ref{2.12}) and (\ref{2.13}), 
there is a positive constant
$C$ such that
\begin{equation} \label{2.14}
\left| \frac{d\lambda_j}{dt} \right| \: \le \: C \: 
\parallel \dot{c}(t) 
\parallel_{c(t)} \: \left( \lambda_j^2 \: + 
K \right)^{1/2}.
\end{equation}
Integration gives equation (\ref{2.8}).
The proposition follows. 
\end{pf}

For some basic
facts about infranilmanifolds, we refer to
\cite[Section 3]{Lott (1999)}.
Let $N$ be a simply-connected connected nilpotent Lie group. Let 
$\Gamma$ be a discrete subgroup of
$\Aff(N)$ which acts freely and cocompactly on $N$. Put $Z \: = \: \Gamma
\backslash N$, an infranilmanifold.
There is a canonical flat linear connection $\nabla^{aff}$ on $TZ$. Put
$\widehat{\Gamma} \: = \: \Gamma \: \cap \: N$, a cocompact subgroup
of $N$. There is a short exact sequence
\begin{equation} \label{2.15}
1 \longrightarrow \widehat{\Gamma} \longrightarrow \Gamma \longrightarrow F
\longrightarrow 1,
\end{equation}
with $F$ a finite group.
Put $\widehat{Z} \: = \: \widehat{\Gamma} \backslash N$,
a nilmanifold which finitely covers $Z$ with covering group $F$.

Let $g^{TZ}$ be a Riemannian metric on $Z$ which is parallel
with respect to
$\nabla^{aff}$. Let us discuss the condition for $Z$ to be spin.
Suppose first that $Z$ is spin. Choose a spin structure on $Z$. Fix the
basepoint $z_0 \: = \: \Gamma \: e \in Z$.  
As $\nabla^{aff}$ preserves $g^{TZ}$, 
its holonomy lies in $\SO(n)$. 
Hence $\nabla^{aff}$ lifts to a flat connection on the principal
$\Spin(n)$-bundle, which we also denote by
$\nabla^{aff}$. There is a corresponding holonomy representation
$\Gamma \rightarrow \Spin(n)$. 

Conversely, suppose that we do not know {\it a priori} if $Z$ is spin. 
Suppose that the affine holonomy $\Gamma \rightarrow F \rightarrow \SO(n)$
lifts to a homomorphism $\Gamma \rightarrow \Spin(n)$. Naturally,
the existence of this lifting is independent of the particular choice of
$g^{TZ}$. Then
there is a corresponding spin structure on $Z$ with principal bundle
$\Gamma \backslash (N \times \Spin(n))$. The different spin structures
on $Z$ correspond to different lifts of $\Gamma \rightarrow \SO(n)$ to
$\Gamma \rightarrow \Spin(n)$. 
These are labelled by
$\HH^1(\Gamma; \Z_2) \: \cong \: \HH^1(Z; \Z_2)$.
Note that there are examples of 
nonspin flat manifolds \cite{Auslander-Szczarba (1962)}.
Also, even if $Z$ is spin and has a fixed spin structure,
the action of $\Aff(Z)$ on $Z$ generally does
not lift to the principal $\Spin(n)$-bundle, as can be seen for the
$\SL(n, \Z)$-action on $Z \: = \: T^n$.

Now let $G$ be either $\SO(n)$ or $\Spin(n)$. Let $V$ be a $G$-Clifford module.
Suppose that $Z$ has a $G$-structure. If $G \: = \: \SO(n)$ then
we have the affine holonomy 
homomorphism $\rho \: : \: \Gamma \rightarrow \SO(n)$. 
If $G \: = \: \Spin(n)$ then
we have a given lift of it to $\rho \: : \: \Gamma \rightarrow
\Spin(n)$. In either case, there is an action of $\Gamma$ on $V$
coming from $\Gamma \stackrel{\rho}{\rightarrow} G \rightarrow
\Aut(V)$.
The vector bundle $E^Z$ can now be written as $E^Z \: = \: \Gamma \backslash
(N \times V)$.
We see that the vector space of sections of $E^Z$ 
which are parallel with respect to
$\nabla^{aff}$ is isomorphic to $V^\Gamma$, the subspace of $V$ which is
fixed by the action of $\Gamma$.

If $V$ is the spinor representation of $G \: = \: \Spin(n)$ then let us
consider the conditions for $V^\Gamma$ to be nonzero.  First, as the
restriction of $\rho \: : \: \Gamma \rightarrow \Spin(n)$ to 
$\widehat{\Gamma}$ maps $\widehat{\Gamma}$ to $\pm 1$, we must have 
$\rho \big|_{\widehat{\Gamma}} \: =  \: 1$. Given this, the homomorphism
$\rho$ factors through a homomorphism $F \rightarrow \Spin(n)$. Then
we have $V^\Gamma \: = \: V^F$. This may be nonzero even if the
homomorphism $F \rightarrow \Spin(n)$ is nontrivial.

Returning to the case of general $V$, 
as $g^{TZ}$ is parallel with respect to $\nabla^{aff}$, the operator $D^Z$
preserves the space $V^\Gamma$ of affine-parallel sections of $E^Z$. Let
$D^{inv}$ be the restriction of $D^Z$ to $V^\Gamma$. 

\begin{proposition} \label{prop2}
There are positive constants $A$, $A^\prime$ and $C$ depending only on
$\dim(Z)$ and $V$ such that if $\parallel R^Z \parallel_\infty \:
\diam(Z)^2 \: \le \: A^\prime$ then the spectrum $\sigma(D^Z)$ of $D^Z$ 
satisfies
\begin{align} \label{2.16}
& \sigma(D^Z) \: \cap \: [ - \: 
\left( A \: \diam(Z)^{-2} \: - \: C \: \parallel R^Z \parallel_\infty
\right)^{1/2},
\left( A \: \diam(Z)^{-2} \: - \: C \: \parallel R^Z \parallel_\infty
\right)^{1/2}
] \: = \notag \\
& \sigma(D^{inv}) 
\: \cap \: [ - \: 
\left( A \: \diam(Z)^{-2} \: - \: C \: \parallel R^Z \parallel_\infty
\right)^{1/2},
\left( A \: \diam(Z)^{-2} \: - \: C \: \parallel R^Z \parallel_\infty
\right)^{1/2}
].
\end{align}
\end{proposition}
\begin{pf}
As $D^Z$ is diagonal with respect to the orthogonal decomposition
\begin{equation} \label{2.17}
C^\infty(Z; E^Z) \: = \: V^\Gamma \: \oplus \: \left( V^\Gamma \right)^\perp,
\end{equation}
it is enough to show that there are constants $A$, $A^\prime$ and $C$ as in
the statement of the proposition such that the eigenvalues of 
$(D^Z)^2 \big|_{\left( V^\Gamma \right)^\perp}$ are greater than
$A \: \diam(Z)^{-2} \: - \: C \: \parallel R^Z \parallel_\infty$.
As in the proof of \cite[Proposition 2]{Lott (1999)}, we can reduce to the
case when $F \: = \: \{e\}$, i.e. $Z$ is a nilmanifold $\Gamma \backslash N$.
Then 
\begin{equation} \label{2.18}
C^\infty(Z; E^Z) \: \cong \: \left( C^\infty(N) \: \otimes \: V 
\right)^\Gamma.
\end{equation}
Using an orthonormal frame $\{e_i\}_{i=1}^{dim(Z)}$ for the Lie algebra
${\frak n}$ as in the proof of \cite[Proposition 2]{Lott (1999)}, we can write
\begin{equation} \label{2.19}
\nabla^{aff}_{e_i} \: = \: e_i \: \otimes \: \Id.
\end{equation}
and
\begin{equation} \label{2.20}
\nabla^Z_{e_i} \: = \: \left( e_i \: \otimes \: \Id.  \right) \: + \:
\left( \Id. \: \otimes \: \frac{1}{2} \: \sum_{a,b=1}^{dim(Z)} 
\omega_{abi} \: \sigma^{ab} \right).
\end{equation}
The rest of the proof now proceeds as in that of 
\cite[Proposition 2]{Lott (1999)}, to which we refer for details.
\end{pf}

\section{Collapsing to a Smooth Base} \label{sect3}

For background information about superconnections and their applications, 
we refer to
\cite{Berline-Getzler-Vergne (1992)}.
Let $M$ be a connected closed oriented Riemannian
manifold which is the total space of a Riemannian
submersion $\pi \: : \: M \rightarrow B$. 
Suppose that $M$ has a $G^M$-structure
and that $V^M$ is a $G^M$-Clifford module, as in Section \ref{sect2}.
If $G^M \: = \: \SO(n)$, put $G^Z \:= \: \SO(\dim(Z))$ and 
$G^B \:= \: \SO(\dim(B))$.
If $G^M \: = \: \Spin(n)$, put $G^Z \:= \: \Spin(\dim(Z))$ and 
$G^B \:= \: \Spin(\dim(B))$. As a fiber $Z_b$ has a trivial normal bundle
in $M$, it admits a $G^Z$-structure. Fixing an orientation of $T_bB$ fixes
the $G^Z$-structure of $Z_b$. Note, however, that $B$ does not
necessarily have a $G^B$-structure.  For example, if $M$ is oriented then
$B$ is not necessarily oriented, as is shown in the example of
$S^1 \times_{\Z_2} S^2 \rightarrow \R P^2$,
where the generator of
$\Z_2$ acts on $S^1$ by complex conjugation and on $S^2$ by the
antipodal map.  And if $M$ is spin then $B$ is not
necessarily spin, as is shown in the example of $S^5 \rightarrow \C P^2$.
What is true is that if the vertical tangent bundle $TZ$, a vector bundle
on $M$, has a $G^Z$-structure then $B$ has a $G^B$-structure.

Put $E^M \: = \: P \times_{G^M} V^M$.
There is a 
Clifford bundle $C$ on $B$ with the property that
$C^\infty(B; C) \: \cong \: C^\infty(M; E^M)$ 
\cite[Section 9.2]{Berline-Getzler-Vergne (1992)}. If $\dim(Z) \: > \: 0$
then $\dim(C) \: = \: \infty$.
To describe $C$ more explicitly,
let $V^M \: = \: \bigoplus_{l \in L} V^B_l \: \otimes \: V^Z_l$ be the
decomposition of $V_M$ into irreducible representations of
$G^B \times G^Z \subset G^M$. \\ \\
{\bf Examples : } \\
1. If $G^M \: = \: \Spin(n)$ and $V^M$ is the spinor representation then
$V^B$ and $V^Z$ are spinor representations. \\
2. If $G^M \: = \: \SO(n)$ and $V^M \: = \: \Lambda^*(\R^n) \otimes_\R \C$
then $V^B \: = \: \Lambda^*(\R^{dim(B)}) \otimes_\R \C$ and
$V^Z \: = \: \Lambda^*(\R^{dim(Z)}) \otimes_\R \C$. \\

Let $U$ be a contractible open subset of $B$. Choose
an orientation on $U$. For $b \in U$,
let $E^Z_{b,l}$ be the vector bundle on $Z_b$ associated to 
the pair $(G^Z, V^Z_l)$.
Then $E^M \big|_{Z_b} \cong \bigoplus_{l \in L} V^B_l \: \otimes \: E^Z_{b,l}$.
The vector bundles $\{E^Z_{b,l}\}_{b \in U}$ are the fiberwise restrictions
of a vector bundle $E^Z_l$ on $\pi^{-1}(U)$, a vertical ``spinor'' bundle.
There is a pushforward vector
bundle $W_l$ on $U$
whose fiber $W_{l,b}$ over $b \in U$ 
is $C^\infty(Z_b; E^Z_{b,l})$.
If $\dim(Z) \: > \: 0$
then $\dim(W_l) \: = \: \infty$. 
There are Hermitian inner products $\{h^{W_l} \}_{l \in L}$
on $\{W_l\}_{l \in L}$ induced from the vertical Riemannian metric $g^{TZ}$. 
Furthermore, there
are Clifford bundles $\{C_l\}_{l \in L}$ on $U$ for which
the fiber $C_{l,b}$ of $C_l$ over $b \in U$ is isomorphic to 
$V^B_l \: \otimes \: W_{l,b}$.
By construction, $C^\infty \left( Z_b; E^M \big|_{Z_b} \right)
\cong \bigoplus_{l \in L} C_{l,b}$. The Clifford bundles
$\{C_l\}_{l \in L}$ exist globally on $B$ and $C \: = \: \bigoplus_{l \in L}
C_l$. 
The Dirac-type operator $D^M$ decomposes as
$D^M \: = \: \bigoplus_{l \in L} D^M_l$, where $D^M_l$ acts on
$C^\infty(B; C_l)$. 

In order to write $D^M_l$ explicitly,
let us recall the Bismut
superconnection on $W_l$. We will deal with each $l \in L$ separately and so
we drop the subscript $l$ for the moment. We use the notation of
\cite[Section III(c)]{Bismut-Lott (1995)} 
to describe the local geometry of the fiber
bundle $M \rightarrow B$, and the Einstein summation convention. Let
$\nabla^{TZ}$ denote the Bismut connection on $TZ$
\cite[Proposition 10.2]{Berline-Getzler-Vergne (1992)}, which we extend to a 
connection on $E^Z_{l}$.
The Bismut superconnection on $W$ 
\cite[Proposition 10.15]{Berline-Getzler-Vergne (1992)} is of the form
\begin{equation} \label{3.1}
A \: = \: D^W \: + \: \nabla^W \: - \: \frac{1}{4} \: c(T).
\end{equation}
Here $D^W$ is the fiberwise Dirac-type operator and has the form
\begin{equation} \label{3.2}
D^W \: = \: - \: i \: \gamma^j \: \nabla^{TZ}_{e_j} \: 
= \: - \: i \: \gamma^j \left( e_j \: + \: \frac{1}{2} \:
\omega_{pqj} \: \sigma^{pq} \right) .
\end{equation}
Next, $\nabla^W$ is a Hermitian connection on $W$ given by
\begin{equation} \label{3.3}
\nabla^W \: = \: \tau^\alpha \left( \nabla^{TZ}_{e_\alpha} 
\: - \: \frac{1}{2} \: \omega_{\alpha j j} \right)
 \: = \: \tau^\alpha \left( e_\alpha 
\: + \: \frac{1}{2} \: \omega_{jk\alpha} \: \sigma^{jk} 
\: - \: \frac{1}{2} \: \omega_{\alpha j j} \right).
\end{equation}
Finally, 
\begin{equation} \label{3.4}
c(T) \: = \: i \: \omega_{\alpha \beta j} \: \gamma^j \:
\tau^\alpha \: \tau^\beta.
\end{equation}

The superconnection $A$ can be ``quantized'' into an operator $D^A$
on $C^\infty(B; V^B \: \otimes \: W)$. 
Explicitly,
\begin{align} \label{3.5}
D^A \: = \: & - \: i \: \gamma^j \left( e_j \: + \: \frac{1}{2} \:
\omega_{pqj} \: \sigma^{pq} \right) \notag \\
& - \: i \: \gamma^\alpha \left( e_\alpha 
\: + \: \frac{1}{2} \: \omega_{\beta \gamma \alpha} \: 
\sigma^{\beta \gamma}
\: + \: \frac{1}{2} \: \omega_{jk\alpha} \: \sigma^{jk} 
\: - \: \frac{1}{2} \: \omega_{\alpha j j} \right) \notag \\
& + \: i \: \frac{1}{2} \: \omega_{\alpha \beta j} \: \gamma^j \:
\sigma^{\alpha \beta}.
\end{align}
Let ${\cal V} \: \in \: \End(C_l)$ be the self-adjoint operator given by
\begin{equation} \label{3.6}
{\cal V} \: = \: - \: i \: \left(
\omega_{\alpha jk} \: \gamma^k \: \sigma^{\alpha j} \: + \: \frac{1}{2} \:
\omega_{\alpha jj} \: \gamma^\alpha \: + \:
\omega_{\alpha \beta j}  \: (\gamma^j \: 
\sigma^{\alpha \beta} \:
 \: + \:
\gamma^\alpha \: \sigma^{j \beta}) \right).
\end{equation}
Then restoring the index $l$ everywhere,
\begin{equation} \label{3.7}
D^M_l \: = \:  D^{A_l} \: + \: {\cal V}_l.
\end{equation}
\noindent
{\bf Examples : } \\
1. If $G^M \: = \: \Spin(n)$ and $V^M$ is the spinor representation then
${\cal V} \: = \: 0$.\\
2. If $G^M \: = \: \SO(n)$ and $V^M \: = \: \Lambda^*(\R^n) \otimes_\R \C$ then
\begin{equation} \label{3.8}
{\cal V} \: = \: - \: \frac{1}{4} \: i \: \left(
\omega_{\alpha jk} \: \gamma^k \: [\widehat{\gamma}^{\alpha},
\widehat{\gamma}^j ] \: + \:
\omega_{\alpha \beta j} (\gamma^j \: 
[\widehat{\gamma}^{\alpha},
\widehat{\gamma}^{\beta} ] \: + \:
\gamma^\alpha \: [\widehat{\gamma}^{j},
\widehat{\gamma}^{\beta} ] ) \right).
\end{equation}

Now suppose that $\pi \: : \: M \rightarrow B$ is a Riemannian affine
fiber bundle. Then $E^M \big|_{Z_b}$ inherits a flat connection
from the flat affine connections on 
$\{E^Z_{b,l}\}_{l \in L}$. Let $E^B$ be the Clifford bundle on $B$ whose fiber
over $b \in B$ is the space of parallel sections
of $E^M \big|_{Z_b}$. Then $D^M$ restricts to a first-order differential
operator $D^B$ on $C^\infty(B; E^B)$. 

Given $b \in U$ and $l \in L$, 
let $W^{inv}_{l,b}$ be the finite-dimensional subspace of 
$W_{l,b}$ consisting of affine-parallel elements of $C^\infty(Z_b; E^Z_{b,l})$.
From the discussion in Section \ref{sect2}, $W^{inv}_{l,b}$ is 
isomorphic to $\left( V^Z_l \right)^\Gamma$. The vector spaces 
$W^{inv}_{l,b}$ fit together to form a 
finite-dimensional subbundle $W^{inv}_l$ of $W_l$.
There is a corresponding finite-dimensional
Clifford subbundle $C^{inv}_l$ of $C_l$ whose fiber over $b \in U$ is
isomorphic to $V^B_l \: \otimes \: W^{inv}_{l,b}$. Again, $C^{inv}_l$ exists
globally on $B$.
Then $E^B \: = \: \bigoplus_{l \in L} C^{inv}_l$.
Let $D^B_l$ be the restriction of $D^M_l$ to $C^\infty(B; C^{inv}_l)$.
Then
\begin{equation} \label{3.9}
D^B \: = \: \bigoplus_{l \in L} D^B_l.
\end{equation}
The superconnection $A_l$ restricts to an superconnection $A^{inv}_l$ on
$W^{inv}_l$, the endomorphism ${\cal V}_l$ restricts to an endomorphism of
$C^{inv}_l$
and $D^M_l$ restricts to the first-order differential operator
\begin{equation} \label{3.10}
D^B_l \:= \: D^{A_l^{inv}} \: + \: {\cal V}_l^{inv}
\end{equation}
on $C^\infty(B; C_l^{inv})$. \\ \\
\noindent
{\bf Proof of Theorem \ref{thm1} :} \\
The operator $D^M_l$ is diagonal with respect to the orthogonal decomposition
\begin{equation} \label{3.11}
C_l \: = \: C^{inv}_l \: \oplus \: \left( C^{inv}_l \right)^\perp.
\end{equation}
Thus it suffices to show that there are constants $A$, $A^\prime$ and $C$ 
such that the spectrum of $\sigma(D^M_l)$, when restricted
to $\left( C^{inv}_l \right)^\perp$, is disjoint from (\ref{1.1}).  

For simplicity, we drop the subscript $l$. Given 
$\eta \in \: C^\infty \left( B; \left( C^{inv} \right)^\perp \right) \subset
C^\infty(M; E^M)$, it is enough to show that for suitable constants,
\begin{equation} \label{3.12}
\langle D^M \eta, D^M \eta \rangle \: \ge \:
\left( \const \: \diam(Z)^{-2} \: - \: \const \: \left( 
\parallel R^M \parallel_\infty \: + \: \parallel \Pi \parallel_\infty^2
\: + \: \parallel T \parallel_\infty^2 \right) \right)
\: \langle \eta, \eta \rangle. 
\end{equation}
Using (\ref{2.4}), it is enough to show that 
\begin{equation} \label{3.13}
\langle \nabla^M \eta, \nabla^M \eta \rangle \: \ge \:
\left( \const \: \diam(Z)^{-2} \: - \: \const \: \left( 
\parallel R^M \parallel_\infty \: + \: 
\parallel \Pi \parallel_\infty^2
\: + \: \parallel T \parallel_\infty^2 \right) \right)
\: \langle \eta, \eta \rangle. 
\end{equation}

We can write $\nabla^M \: = \: \nabla^V \: + \: \nabla^H$, where
\begin{equation} \label{3.14}
\nabla^V \: : \: C^\infty(M; E^M) \rightarrow C^\infty(M; T^* Z \otimes E^M) 
\end{equation}
denotes covariant differentiation in the vertical direction and
\begin{equation} \label{3.15}
\nabla^H \: : \: C^\infty(M; E^M) 
\rightarrow C^\infty(M; \pi^* T^*B \otimes E^M) 
\end{equation}
denotes covariant differentiation in the horizontal direction.
Then
\begin{align} \label{3.16}
\langle \nabla^M \eta, \nabla^M \eta \rangle \: & = \:
\langle \nabla^V \eta, \nabla^V \eta \rangle \: + \:
\langle \nabla^H \eta, \nabla^H \eta \rangle \notag \\
& \ge \: \langle \nabla^V \eta, \nabla^V \eta \rangle  \notag \\
& = \: \int_B \int_{Z_b} \big| \nabla^V \eta \big|^2 (z) \: d\vol_{Z_b} \:
d\vol_B. 
\end{align}
On a given fiber $Z_b$, we have 
\begin{equation} \label{3.17}
E^M \big|_{Z_b} \: \cong \: V^B \: \otimes \: E^Z_b.
\end{equation}
Hence we can also use the Bismut connection $\nabla^{TZ}$ to vertically
differentiate sections of $E^M$. That is, we can define
\begin{equation} \label{3.18}
\nabla^{TZ} \: : \: C^\infty(M; E^M) \rightarrow C^\infty(M; T^*Z \otimes E^M).
\end{equation}
Explicitly, with respect to a local framing,
\begin{equation} \label{3.19}
\nabla^{TZ}_{e_j} \: = \: e_j \: \eta \: + \: \frac{1}{2} \:
\omega_{pqj} \: \sigma^{pq} \:
\eta
\end{equation}
and
\begin{equation} \label{3.20}
\nabla^V_{e_j} \: = \: e_j \: \eta \: + \: \frac{1}{2} \:
\omega_{pqj} \: \sigma^{pq} \: \eta \: + \:
\omega_{\alpha kj} \: \sigma^{\alpha k} \: \eta \: + \: \frac{1}{2} \:
\omega_{\alpha \beta j} \: \sigma^{\alpha \beta} \: \eta.
\end{equation}
Then from (\ref{3.16}), (\ref{3.19}) and (\ref{3.20}),
\begin{equation} \label{3.21}
\langle \nabla^M \eta, \nabla^M \eta \rangle \: \ge \:
\int_B \left[ \int_{Z_b} \big| \nabla^{TZ} \eta \big|^2 (z) \: - \:
\const \left( \parallel T_b \parallel^2 \: + \: 
\parallel \Pi_b \parallel^2 \right) \big| \eta(z) \big|^2 \right] \:
d\vol_{Z_b} \: d\vol_B.
\end{equation}

Thus it suffices to bound $\int_{Z_b} \big| \nabla^{TZ} \eta \big|^2 (z)
\: d\vol_{Z_b}$ from below on a given fiber $Z_b$ in terms of
$\langle \eta, \eta \rangle_{Z_b}$, under the assumption that
$\eta \in (W^{inv}_b)^\perp$. Using the Gauss-Codazzi equation, we can
estimate $\parallel R^{Z_b} \parallel_\infty$ in terms of
$\parallel R^{M} \parallel_\infty$ and 
$\parallel \Pi \parallel_\infty^2$. Then the desired bound on
$\int_{Z_b} \big| \nabla^{TZ} \eta \big|^2 (z)
\: d\vol_{Z_b}$ follows from Proposition \ref{prop2}. $\square$ \\ \\
{\bf Proof of Theorem \ref{thm2} :} \\
Let $g_0^{TM}$ denote the Riemannian metric on $M$.
From Proposition \ref{prop1},
if a Riemannian metric $g_1^{TM}$ on $M$ is close
to $g_0^{TM}$ in $({\cal M}(M, 2K), d)$
then applying the function
$x \rightarrow \sinh^{-1} \left( \frac{x}{\sqrt{2K}} \right)$ to
$\sigma(D^M, g_0^{TM})$ gives a collection of numbers which is close to
that obtained by applying
$x \rightarrow \sinh^{-1} \left( \frac{x}{\sqrt{2K}} \right)$ to
$\sigma(D^M, g_1^{TM})$.
We will use the geometric results of \cite{Cheeger-Fukaya-Gromov (1992)} to
find a metric $g_2^{TM}$ on $M$ which is close to $g_0^{TM}$
and to which we can apply Theorem \ref{thm1}.

First, as in \cite[(2.4.1)]{Cheeger-Fukaya-Gromov (1992)}, by the
smoothing results of Abresch and others
\cite[Theorem 1.12]{Cheeger-Fukaya-Gromov (1992)},
for any $\epsilon \: > \: 0$ we can find metrics
on $M$ and $B$ which are $\epsilon$-close in the $C^1$-topology
to the original metrics such that
the new metrics satisfy $\parallel \nabla^i R \parallel_\infty \: \le \:
A_i(n, \epsilon)$ for some appropriate sequence 
$\{A_i(n, \epsilon)\}_{i=0}^\infty$.
Let $g_1^{TM}$ denote the new metric on $M$.
In the proof of the smoothing result, such as using the Ricci flow
\cite[Proposition 2.5]{Rong (1996)}, 
one obtains an explicit smooth $1$-parameter family of metrics on $M$ in
${\cal M}(M, K^\prime)$, for some $K^\prime \: > \: K$, going from
$g_0^{TM}$ to $g_1^{TM}$. We can approximate this family by a
piecewise-analytic family. Hence one obtains an upper bound on 
$d \left( g_0^{TM},
g_1^{TM} \right)$ in ${\cal M}(M, K^\prime)$,
for some $K^\prime \: > \: K$, which depends on $K$ and
is proportionate to $\epsilon$. (Note that $d$ is essentially the same as
the $C^0$-metric on ${\cal M}(M, K^\prime)$.)
By rescaling, we may assume that
$\parallel R^M \parallel_\infty \: \le \: 1$,
$\parallel R^B \parallel_\infty \: \le \: 1$ and $\inj(B) \ge 1$.
We now
apply \cite[Theorem 2.6]{Cheeger-Fukaya-Gromov (1992)}, with $B$ fixed.
It implies that
there are positive constants $\lambda(n)$ and $c(n,\epsilon)$ so that if
$d_{GH}(M, B) \: \le \: \lambda(n)$ then there is a fibration
$f : M \rightarrow B$ such that\\
1. $\diam\left( f^{-1}(b) \right) \: \le \: c(n, \epsilon) \: d_{GH}(M, B)$.\\
2. $f$ is a $c(n, \epsilon)$-almost Riemannian submersion.\\
3. $\parallel \Pi_{f^{-1}(b)} \parallel_\infty \: \le \: 
c(n, \epsilon)$.\\
As in \cite{Fukaya (1989)}, the Gauss-Codazzi equation, 
the curvature bound on $M$ and the
second fundamental form bound on $f^{-1}(b)$ imply a
uniform bound on 
$\left\{\parallel R^{f^{-1}(b)} \parallel_\infty \right\}_{b \in B}$. 
Along with
the diameter bound on $f^{-1}(b)$, this implies that if $d_{GH}(M, B)$ is
sufficiently small then $f^{-1}(b)$ is almost flat. 

From
\cite[Propositions 3.6 and 4.9]{Cheeger-Fukaya-Gromov (1992)}, we can find
another metric $g_2^{TM}$ on $M$ which is $\epsilon$-close to $g_1^{TM}$ 
in the $C^1$-topology so that the fibration
$f : M \rightarrow B$ gives $M$ the structure of a Riemannian 
affine fiber bundle.
Furthermore, by \cite[Proposition 4.9]{Cheeger-Fukaya-Gromov (1992)},
there is a sequence $\{A^\prime_i(n, \epsilon)\}_{i=0}^\infty$ so that we may
assume that $g_1^{TM}$ and $g_2^{TM}$ are close in the sense that
\begin{equation} \label{3.22}
\parallel \nabla^i \left( g_1^{TM} - g_2^{TM} \right) \parallel_\infty
\: \le \: A^\prime_i(n, \epsilon) \: d_{GH}(M,B),
\end{equation}
where the covariant derivative in (\ref{3.22})
is that of the Levi-Civita connection of
$g_2^{TM}$. 
Then we can interpolate linearly between $g_1^{TM}$ and $g_2^{TM}$ within
${\cal M}(M, K^{\prime \prime})$ for some $K^{\prime \prime} > K^\prime$,
and obtain an upper bound on $d \left( g_1^{TM}, g_2^{TM} \right)$ in 
${\cal M}(M, K^{\prime \prime})$ which is proportionate to $\epsilon$.
From \cite[Theorem 2.1]{Rong (1996)}, we can take $K^{\prime \prime} \: = \:
2 K$ (or any number greater than $K$).

We now apply Theorem \ref{thm1} to the Riemannian 
affine fiber bundle with metric
$g_2^{TM}$. It remains
to estimate the geometric terms appearing in (\ref{1.1}). 
We have an estimate
on $\parallel \Pi \parallel_\infty$ as above. Applying
O'Neill's formula \cite[(9.29)]{Besse (1987)} 
to the Riemannian affine fiber bundle, 
we can estimate $\parallel T \parallel_\infty^2$ in terms
of $\parallel R^M \parallel_\infty$ and $\parallel R^B \parallel_\infty$.
Putting this together, the theorem follows. $\square$

\section{Collapsing to a Singular Base} \label{sect4}

Let ${\frak p} : P \rightarrow M$ be the principal $G$-bundle
of Section \ref{sect2}.
Let $\{{\frak Y}_j\}_{j=1}^n$ be the horizontal vector fields on $P$ such
that $\theta({\frak Y}_j) = e_j$.
Put
$D^P \: = \: - \: i \: \sum_{j=1}^n \gamma^j \: {\frak Y}_j$, acting on
$C^\infty(P) \otimes V$.

There is an isomorphism 
$C^\infty(M; E^M) \: \cong \: \left( C^\infty(P) \otimes V \right)^G$.
Under this isomorphism, 
$D^M \: \cong \: D^P \big|_{\left( C^\infty(P) \otimes V \right)^G}$.
The Bochner-type equation (\ref{2.4}) becomes
\begin{equation} \label{4.1}
(D^M)^2 \: \cong \: - \: \sum_{j=1}^n {\frak Y}_j^2 
\: - \: \frac{1}{8} \: \sum_{a,b,i,j=1}^n
({\frak p}^* R^M)_{abij} \: 
(\gamma^i \: \gamma^j \: - \: \gamma^j \: \gamma^i) \:
\sigma^{ab} 
\end{equation}
when acting on $\left( C^\infty(P) \otimes V \right)^G$.

Let $\{x_a\}_{a=1}^{\dim(G)}$ be a basis for the Lie algebra ${\frak g}$ which
is orthonormal with respect to the negative of the Killing form.
Let $\{{\frak Y}_a\}_{a=1}^{\dim(G)}$ be the corresponding vector fields on 
$P$. Then $- \sum_{a=1}^{\dim(G)} {\frak Y}_a^2$ acts on 
$\left( C^\infty(P) \otimes V \right)^G$ as $c_V \in \R$, the Casimir of
the $G$-module $V$. 
Give $P$ the Riemannian metric $g^{TP}$ with the property that
$\{{\frak Y}_j, {\frak Y}_a\}$ forms an orthonormal basis of vector fields. Let
$\triangle^P$ denote the corresponding (nonnegative)
scalar Laplacian on $P$, extended to act
on $C^\infty(P) \otimes V$. Then when acting on 
$\left( C^\infty(P) \otimes V \right)^G$, equation (\ref{4.1}) is equivalent to
\begin{equation} \label{4.2}
(D^M)^2 \: \cong \: \triangle^P \: - \: \frac{1}{8} \: \sum_{a,b,i,j=1}^n
({\frak p}^* R^M)_{abij} \: 
(\gamma^i \: \gamma^j \: - \: \gamma^j \: \gamma^i) \:
\sigma^{ab} \: - \: c_V \: \Id.
\end{equation}

\begin{definition} \label{def3}
A $G$-equivariant Riemannian affine fiber bundle structure on $P$
consists of a Riemannian affine fiber bundle structure
$\check{\pi} : P \rightarrow \check{X}$ which is $G$-equivariant.
\end{definition}

Given a $G$-equivariant Riemannian affine fiber bundle, let
$\check{Z}$ be the fiber of
$\check{\pi} : P \rightarrow \check{X}$, an infranilmanifold.
For collapsing purposes 
it suffices to take $\check{Z}$ to be a nilmanifold $\Gamma \backslash N$
\cite[(7.2)]{Cheeger-Fukaya-Gromov (1992)}. 
We assume hereafter that this is the case.
Put $X \: = \: \check{X}/G$, a possibly singular space.
As $N$ acts isometrically in a neighborhood of a given fiber $\check{Z}$
and preserves the horizontal subspaces of $P \rightarrow M$,
it follows that 
the vector fields $\{{\frak Y}_j\}_{j=1}^n$ are projectable with respect to
$\check{\pi}$
and push forward to vector fields $\{{\cal X}_j\}_{j=1}^n$ on $\check{X}$.
Put $D^{\check{X}} \: = \: - \: i \:
\sum_{j=1}^n \gamma^j \: {\cal X}_j$, acting on
$C^\infty(\check{X}) \otimes V$. Let $v \in C^\infty(\check{X})$ be given
by $v(\check{x}) \: = \: \vol(\check{Z}_{\check{x}})$. We give
$C^\infty(\check{X}) \otimes V$ the weighted $L^2$-inner product with
respect to the weight function $v$. 

We recall that there is a notion of a pseudodifferential operator being 
transversally elliptic with respect to the action of a Lie group $G$
\cite[Definition 1.3]{Atiyah (1974)}.

\begin{lemma} \label{lemma1}
$D^{\check{X}}$ is transversally elliptic on $\check{X}$. 
\end{lemma}
\begin{pf}
Let $s(D^{\check{X}}) \in C^\infty(T^* \check{X})\: \otimes \: \End(V)$ 
denote the symbol of $D^{\check{X}}$.
Suppose that $\xi \in T^*_{\check{x}} \check{X}$ satisfies
$\xi(\check{v}) \: = \: 0$ for all
$\check{v} \in T_{\check{x}} \check{X}$ which lie in the image of the 
representation of ${\frak g}$ by vector fields on $\check{X}$.
Then if $p \in \check{\pi}^{-1}(\check{x})$, we have that
$(\check{\pi}^* \xi)(r) \: = \: 0$ for all
$r \in T_p P$ which lie in the image of the 
representation of ${\frak g}$ by vector fields on $P$. In other words,
$\check{\pi}^* \xi$ is horizontal.
Suppose in addition that $s(D^{\check{X}})(\xi) \: = \: 0$. Then
$s(D^P)(\check{\pi}^* \xi) \: = \: 0$. As $D^P$ is horizontally elliptic,
it follows that $\check{\pi}^* \xi \: = \: 0$. Thus $\xi \: = \: 0$, which
proves the lemma.
\end{pf}

\begin{definition} \label{def4}
For notation, write $C^\infty(X; E^X) \: = \: 
\left( C^\infty(\check{X}) \otimes V \right)^G$.
Let $D^X$ be the restriction of $D^{\check{X}}$ to 
$C^\infty(X; E^X)$.
\end{definition}

It will follow from the proof of the
next theorem that $D^X$ is self-adjoint on the Hilbert space completion of
$C^\infty(X; E^X)$ with respect to the (weighted)
inner product.  As $D^{\check{X}}$ is transversally 
elliptic, it follows that $D^X$ has a discrete spectrum
\cite[Proof of Theorem 2.2]{Atiyah (1974)}.

Let $\check{\Pi}$ denote the second fundamental forms of the fibers
$\{\check{Z}_{\check{x}} \}_{\check{x} \in \check{X}}$. Let
$\check{T} \in \Omega^2(P; T \check{Z})$ be the 
curvature of the horizontal distribution on the affine fiber bundle
$P \rightarrow \check{X}$.

\begin{theorem} \label{thm6}
There are positive constants $A$, $A^\prime$ and $C$ which only depend on
$n$ and $V$ such that if $\parallel R^{\check{Z}} \parallel_\infty \:
\diam(\check{Z})^2 \: \le \: A^\prime$ then the intersection of 
$\sigma(D^M)$ with
\begin{align} \label{4.3}
[ 
- \left( \: A \: \diam(\check{Z})^{-2} \: - \: \right.
& 
\left. C \: \left( 1 \: + \:
\parallel R^M \parallel_\infty \: + \: {\parallel \check{\Pi} 
\parallel_\infty}^2
\: + \: {\parallel \check{T} 
\parallel_\infty}^2 \right) \right)^{1/2},  \notag  \\
&  
\left( A \: \diam(\check{Z})^{-2} \: - \: C \: \left( 1 \: + \:
\parallel R^M \parallel_\infty \: + \: {\parallel \check{\Pi} 
\parallel_\infty}^2
\: + \: {\parallel \check{T} \parallel_\infty}^2 \right) \right)^{1/2}
] 
\end{align}
equals the intersection of $\sigma(D^X)$ with (\ref{4.3}).
\end{theorem}
\begin{pf}
Let us write
\begin{equation} \label{4.4}
C^\infty(P) \otimes V \: = \: \left( C^\infty(\check{X}) \otimes V \right) 
\oplus \left( C^\infty(\check{X}) \otimes V \right)^\perp,
\end{equation}
where we think of $C^\infty ( \check{X}) \otimes V$ as the elements
of $C^\infty(P) \otimes V$ which are constant along the fibers of the fiber 
bundle $\check{\pi} \: : \: P \rightarrow \check{X}$.
Taking $G$-invariant subspaces, we have an orthogonal decomposition
\begin{equation} \label{4.5}
C^\infty(M; E^M) \: = \: C^\infty(X; E^X) \oplus \left( C^\infty(X; E^X)
\right)^\perp,
\end{equation}
with respect to which $D^M$ decomposes as
\begin{equation} \label{4.6}
D^M \: = \: D^X \oplus D^M 
\big|_{\left( C^\infty(X; E^X)
\right)^\perp}.
\end{equation}
As in the proof of Theorem \ref{thm1},
it suffices to obtain a lower bound on
the spectrum of $(D^M)^2 
\big|_{\left( C^\infty(X; E^X)
\right)^\perp}$.
As
$\left( C^\infty(X; E^X)
\right)^\perp \subset \left( C^\infty(\check{X}) \otimes V \right)^\perp$,
using (\ref{4.2}) it suffices to obtain a lower bound on
the spectrum of 
$\triangle^P \big|_{\left( C^\infty(\check{X}) \otimes V \right)^\perp}$.
This follows from the arguments of the proof of Theorem \ref{thm1}, using the
fact that $\parallel R^P \parallel_\infty \: \le \: \const \left( 1 \: + \:
\parallel R^M \parallel_\infty \right)$. We omit the 
details. In fact, it is somewhat easier than the proof of Theorem \ref{thm1},
since we are now only dealing with the
scalar Laplacian and so can replace Proposition \ref{prop2} by standard
eigenvalue estimates (which just involve a lower Ricci curvature bound); see
\cite{Berard (1988)} and references therein.
\end{pf}
\noindent
{\bf Proof of Theorem \ref{thm3} :} \\
Everything in the proof will be done in a $G$-equivariant way, so we may
omit to mention this explicitly. Let $P_i$ be the principal $G$-bundle of
$M_i$, equipped with a Riemannian metric as in the beginning of the section.
From the $G$-equivariant version of
Gromov's compactness theorem, we obtain a subsequence 
$\{P_i\}_{i=1}^\infty$ which converges in the 
equivariant Gromov-Hausdorff topology
to a $G$-Riemannian manifold $\left( \check{X}, g^{T\check{X}} \right)$ with
a $C^{1,\alpha}$-regular metric. As in \cite[Section 3]{Fukaya (1987)}, the
measure $\chi \: d\vol_{\check{X}}$ 
is a weak-$*$ limit point of the pushforwards of the normalized 
Riemannian measures on $\{P_i\}_{i=1}^\infty$. 
As in \cite[p. 535]{Fukaya (1987)},
after smoothing we may assume that we have $G$-equivariant 
Riemannian affine fiber bundles $\check{\pi}_i : P_i^\prime \rightarrow
\check{X}_i$, with $G$ acting freely on $P_i^\prime$, along with
$G$-diffeomorphisms $\check{\phi_i} \: : \: 
P_i \rightarrow P_i^\prime$ and $\Phi_i \: : \: 
\check{X} \rightarrow \check{X}_i$. Put $M^\prime_i \: = \: P_i^\prime/G$.
Then $\check{\phi_i}$ descends to a diffeomorphism  
$\phi_i \: : \: 
M_i \rightarrow M_i^\prime$ and
we may also assume, as in the proof of Theorem \ref{thm2}, that \\
1. $\phi_i^* g^{TM_i^\prime} \: \in \: {\cal M}(M_i, \const \: K)$, \\
2. $d(\phi_i^* g^{TM_i^\prime}, g^{TM_i}) \: \le \: 
2^{- \: i}$ in ${\cal M}(M_i, \const \: K)$ and \\
3. $\lim_{i \rightarrow \infty}
\Phi_i^* g^{T\check{X}_i} \: = \: g^{T\check{X}}$ in the 
$C^{1,\alpha}$-topology. \\
Using Proposition \ref{prop1}, we can effectively replace $M_i$ by 
$M_i^\prime$ for the purposes of the argument. For simplicity, we relabel
$M_i^\prime$ as $M_i$ and
$P_i^\prime$ as $P_i$.
For the purposes of the limiting argument, 
using Theorem \ref{thm6} and (\ref{4.2}), 
we may replace the spectrum of $|D^{M_i}|$ by the spectrum of
the operator $|D^{X_i}| \: \equiv \:
\sqrt{\triangle^{\check{X}_i} \: + \: {\cal V}_i}$
acting on $C^\infty(X_i, E^{X_i}) \: = \: \left( C^\infty(\check{X}_i) \otimes
V \right)^G$, where
${\cal V}_i$ is the restriction of
\begin{equation} \label{4.7}
- \: \frac{1}{8} \: \sum_{a,b,i,j=1}^n
({\frak p}^* R^{M_i})_{abij} \: 
(\gamma^i \: \gamma^j \: - \: \gamma^j \: \gamma^i) \:
\sigma^{ab} \: - \: c_V \: \Id.
\end{equation}
to the elements of $\left( C^\infty(P_i) \otimes V \right)^G$ which
are constant along the fibers of $\check{\pi}_i : P_i \rightarrow
\check{X}_i$, i.e. to $C^\infty(X_i, E^{X_i})$. 

From the curvature bound, we have a uniform bound on
$\{\parallel {\cal V}_i \parallel_\infty \}_{i=1}^\infty$. Using the
weak-$*$ compactness of the unit ball,
let ${\cal V}$ be a weak-$*$ limit point of 
$\{\Phi_i^* {\cal V}_i \}_{i=1}^\infty$ in
$L^\infty(\check{X}) 
\: \otimes \: \End(V) \: = \: \left( L^1(\check{X}) 
\: \otimes \: \End(V) \right)^*$.
We claim that with this choice of $\check{X}$, $\chi$ and ${\cal V}$, equation
(\ref{1.4}) holds. 

To see this, we use the minimax characterization of
eigenvalues as in \cite[Section 5]{Fukaya (1987)}. Using the diffeomorphisms
$\{\Phi_i\}_{i=1}^\infty$, we identify each $\check{X}_i$ with 
$\check{X}$. We denote by 
$\langle \cdot, \cdot \rangle_{X_i}$ an $L^2$-inner product
constructed using $\Phi_i^* g^{T \check{X}_i}$ and the weight function
$(\check{\pi}_i)_*(d\vol_{P_i})/\int_{\check{X}_i} 
(\check{\pi}_i)_*(d\vol_{P_i})$.
We denote by 
$\langle \cdot, \cdot \rangle_{X}$ an $L^2$-inner product
constructed using $g^{T \check{X}}$ and the weight function
$\chi \: d\vol_{\check{X}}$.
As $\triangle^{\check{X}}$ has a compact resolvent,
it follows that $|D^X|^2$ has a compact resolvent. Then
\begin{equation} \label{4.8}
\lambda_k(|D^X|)^2 \: = \: \inf_W \sup_{\psi \in W-0} 
\frac{\langle d \psi, d \psi \rangle_X \: + \: \langle \psi,
{\cal V} \psi \rangle_X}{\langle \psi, \psi \rangle_X},
\end{equation}
where $W$ ranges over the $k$-dimensional subspaces of the Sobolev space
$H^1(X; E^X)$. Given $\epsilon \: > \: 0$, let $W_\infty$ be
a $k$-dimensional subspace such that
\begin{equation} \label{4.9}
\sup_{\psi \in W_\infty -0} 
\frac{\langle d \psi, d \psi \rangle_X \: + \: \langle \psi,
{\cal V} \psi \rangle_X}{\langle \psi, \psi \rangle_X} \: \le \:
\lambda_k(|D^X|)^2 \: + \epsilon.
\end{equation}
As $\psi \otimes \psi^*$ lies in the finite-dimensional subspace
$W_\infty \otimes W_\infty^*$ of $L^1(\check{X}) 
\otimes \End(V)$, it follows that 
\begin{equation} \label{4.10}
\lim_{i \rightarrow \infty} \langle \psi, {\cal V}_i \psi \rangle_X
\: = \:  
\langle \psi, {\cal V} \psi \rangle_X
\end{equation}
uniformly on 
$\{ \psi \in W_\infty \: : \: \langle \psi, \psi \rangle_X \: = \: 1 \}$.
Then 
\begin{equation} \label{4.11}
\lim_{i \rightarrow \infty} \sup_{\psi \in W_\infty -0} 
\frac{\langle d \psi, d \psi \rangle_{X_i} \: + \: \langle \psi,
{\cal V}_i \psi \rangle_{X_i}}{\langle \psi, \psi \rangle_{X_i}} \: = \:
\sup_{\psi \in W_\infty -0} 
\frac{\langle d \psi, d \psi \rangle_X \: + \: \langle \psi,
{\cal V} \psi \rangle_X}{\langle \psi, \psi \rangle_X}.
\end{equation}
As
\begin{equation} \label{4.12}
\lambda_k(|D^{X_i}|)^2 \: = \: \inf_W \sup_{\psi \in W-0} 
\frac{\langle d \psi, d \psi \rangle_{X_i} \: + \: \langle \psi,
{\cal V}_i \psi \rangle_{X_i}}{\langle \psi, \psi \rangle_{X_i}},
\end{equation}
it follows that
\begin{equation} \label{4.13}
\limsup_{i \rightarrow \infty} \lambda_k(|D^{X_i}|) \: \le \:
\lambda_k(|D^{X}|). 
\end{equation}

We now show that
\begin{equation} \label{4.14}
\liminf_{i \rightarrow \infty} \lambda_k(|D^{X_i}|) \: \ge \:
\lambda_k(|D^{X}|). 
\end{equation}

Along with (\ref{4.13}), this will prove the theorem.
Suppose that (\ref{4.14}) is not true.  Then there is some $\epsilon \: > \: 0$
and some infinite subsequence of $\{M_i\}_{i=1}^\infty$, which we relabel
as $\{M_i\}_{i=1}^\infty$, such that for all $i \in \Z^+$,
\begin{equation} \label{4.15}
\lambda_k(|D^{X_i}|)^2 \: \le \: \lambda_k(|D^{X}|)^2 \: - \: 2 \: \epsilon.
\end{equation}
For each $i \in \Z^+$, let $W_i$ be a $k$-dimensional subspace of
$H^1(X; E^X)$ such that
\begin{equation} \label{4.16}
\sup_{\psi \in W_i -0} 
\frac{\langle d \psi, d \psi \rangle_{X_i} \: + \: \langle \psi,
{\cal V}_i \psi \rangle_{X_i}}{\langle \psi, \psi \rangle_{X_i}} 
\: \le \: \lambda_k(|D^{X_i}|)^2 \: + \: \epsilon.
\end{equation}
Let $\{f_{i,j}\}_{j=1}^k$ be a basis for $W_i$ which is orthonormal
with respect to $\langle \cdot, \cdot \rangle_X$. Then for a given
$j$,
the sequence $\{f_{i,j}\}_{i=1}^\infty$ is bounded in
$H^1(X; E^X)$. After taking a subsequence, which we relabel as
$\{f_{i,j}\}_{i=1}^\infty$,
we can assume that $\{f_{i,j}\}_{i=1}^\infty$ converges weakly in
$H^1(X; E^X)$ to
some $f_{\infty, j}$. Doing this successively for 
$j \in \{1, \ldots, k\}$, we can assume that for each $j$,
$\lim_{i \rightarrow \infty} f_{i,j} \: = \: f_{\infty, j}$ weakly in
$H^1(X; E^X)$. Then from the compactness of the embedding
$H^1(X; E^X) \rightarrow L^2(X; E^X)$, 
we have strong convergence in $L^2(X; E^X)$. In particular,
$\{f_{\infty,j}\}_{j=1}^k$ are orthonormal. Put $W_\infty \: = \:
\spann (f_{\infty,1}, \ldots, f_{\infty,k})$. 

If $w_\infty \: = \: \sum_{j=1}^k c_j \: f_{\infty, j}$ is a nonzero element of
$W_\infty$, put $w_i \: = \: \sum_{j=1}^k c_j \: f_{i, j}$. Then
$\{w_i\}_{i=1}^\infty$ converges weakly to $w_\infty$ in $H^1(X; E^X)$ and
hence converges strongly to $w_\infty$ in $L^2(X; E^X)$.
From a general result about weak limits, we have
\begin{equation} \label{4.17}
\langle w_\infty , w_\infty \rangle_{H^1} \: \le \: 
\limsup_{i \rightarrow \infty} \:
\langle w_i , w_i \rangle_{H^1}.
\end{equation}
Along with the $L^2$-convergence of
$\{w_i\}_{i=1}^\infty$ to $w_\infty$, this implies that
\begin{equation} \label{4.18}
\langle dw_\infty, dw_\infty \rangle_X \: \le \: 
\limsup_{i \rightarrow \infty} \:
\langle dw_i, dw_i \rangle_{X_i}.
\end{equation}
As $w_i \: \otimes \: w_i^*$ converges in $L^1(\check{X}) \: \otimes \:
\End(E)$ to $w_\infty \: \otimes \: w_\infty^*$, we have
\begin{align} \label{4.19}
\lim_{i \rightarrow \infty} \langle w_i, {\cal V}_i w_i \rangle_{X} \: & = \:
\lim_{i \rightarrow \infty} 
\left( \langle w_\infty, {\cal V}_i w_\infty \rangle_{X} \: + \: 
\left( \langle w_i, {\cal V}_i w_i \rangle_{X} \: - \:
\langle w_\infty, {\cal V}_i w_\infty \rangle_{X}\right) \right) \notag \\
& = \:
\langle w_\infty, {\cal V} w_\infty \rangle_{X}.
\end{align}
Then
\begin{equation} \label{4.20}
\sup_{\psi \in W_\infty - 0} 
\frac{\langle d \psi, d \psi \rangle_{X} \: + \: \langle \psi,
{\cal V} \psi \rangle_{X}}{\langle \psi, \psi \rangle_{X}}
\: \le \: \limsup_{i \rightarrow \infty}
\sup_{\psi \in W_i - 0} 
\frac{\langle d \psi, d \psi \rangle_{X_i} \: + \: \langle \psi,
{\cal V}_i \psi \rangle_{X_i}}{\langle \psi, \psi \rangle_{X_i}}.
\end{equation}
Thus from (\ref{4.15}), (\ref{4.16}) and (\ref{4.20}),
\begin{equation} \label{4.21}
\inf_W
\sup_{\psi \in W - 0} 
\frac{\langle d \psi, d \psi \rangle_{X} \: + \: \langle \psi,
{\cal V} \psi \rangle_{X}}{\langle \psi, \psi \rangle_{X}} \: \le \:
\lambda_k(|D^{X}|)^2 \: - \: \epsilon,
\end{equation}
which is a contradiction. This proves the theorem.
$\square$ \\ \\
{\bf Proof of Theorem \ref{thm4} :} \\
Let $\{g_i^{TM}\}_{i=1}^\infty$ be a sequence of Riemannian metrics on $M$
as in the statement of the theorem,
with respect to which $\lambda_k(|D^M|)$ goes to infinity.
Let $P$ be the principal $G$-bundle of $M$ and let $\check{X}$ be the
limit space of Theorem \ref{thm3}, 
a smooth manifold with a $C^{1, \alpha}$-regular
metric.
As the limit space $X \: = \: \check{X}/G$ has diameter $1$, it has positive
dimension.
As in the proof of Theorem \ref{thm3}, after slightly smoothing the metric on
$\check{X}$, there is a $G$-equivariant Riemannian
affine fiber bundle $\check{\pi} \: : \: P \rightarrow \check{X}$ whose
fiber is a nilmanifold $\check{Z}$.
Let $\check{x}$ be a point in a principal
orbit for the $G$-action on $\check{X}$,
with isotropy group $H \subset G$. Then $H$ acts affinely
on the nilmanifold fiber $\check{Z}_{\check{x}}$. In particular, $H$ is
virtually abelian.  The quotient $Z \: = \: \check{Z}_{\check{x}}/H$ is the
generic fiber of the possibly-singular fiber bundle $\pi \: : \:
M \rightarrow X$,
the $G$-quotient of $\check{\pi} \: : \: P \rightarrow \check{X}$.
Then $E^M \big|_Z \: = \: \check{Z}_{\check{x}} \times_H V$. In particular,
the vector space of 
affine-parallel sections of $E^M \big|_Z$ is isomorphic to 
$V^H$. On the other hand, if $C^\infty(X; E^X) \: \neq \: 0$ then
$|D^X|$ has an infinite discrete spectrum.
Theorem \ref{thm3} now implies that $C^\infty(X; E^X) \cong
\left( C^\infty(\check{X}) \otimes V \right)^G$ must be the zero space.
As the orbit $\check{x} \cdot G$ has
a neighborhood consisting of principal orbits, the restriction map
from $\left( C^\infty(\check{X}) \otimes V \right)^G$ to
$\left( C^\infty(\check{x} \cdot G) \otimes V \right)^G$ is surjective. 
However, $\left( C^\infty(\check{x} \cdot G) \otimes V \right)^G$ is
isomorphic to $V^H$. Thus $V^H \: = \: 0$.
This proves the theorem. $\square$

\section{Proof of Theorem \ref{thm5}}

As the proof of Theorem \ref{thm5} is similar to 
\cite[Pf. of Theorem 2]{Lott2 (1999)}, we only indicate the structure of the
proof and the necessary modifications to 
\cite[Pf. of Theorem 2]{Lott2 (1999)}.

The closure $\overline{U_I}$ of an appropriate neighborhood of an end 
has the (affine) structure of
an affine fiber bundle over $[0, \infty)$ with fiber $Z_I$. The vector
bundle $E^B_I$ is the trivial vector bundle over $[0, \infty)$ whose
fiber over $s \in [0, \infty)$ consists of the affine-parallel sections
of $E^M \big|_{\{s\} \times Z_I}$. As in \cite[Section 4]{Lott2 (1999)},
if $U_I$ is sufficiently far out the end then
we can use Propositions \ref{prop1} and \ref{prop2} of the present paper to
construct an embedding of $C^\infty([0, \infty); E^B_I)$ into
$C^\infty \left( \overline{U_I}; E^M 
\big|_{\overline{U_I}} \right)$ whose image
consists of elements with ``bounded energy'' fiberwise restrictions. Let
$P_0$ be the Hilbert space extension of orthogonal projection from 
$\bigoplus_{I=1}^N 
C^\infty \left( \overline{U_I}; E^M \big|_{\overline{U_I}} \right)$ to
$\bigoplus_{I=1}^N C^\infty([0, \infty); E^B_I)$.
By standard arguments as in \cite[Pf. of Proposition 2.1]{Donnelly-Li (1979)},
the essential spectrum of $D^M$ equals that of $D^M_{end}$.
With respect to the decomposition of the Hilbert space into $\Image(P_0) \oplus
\Image(I \:- \: P_0)$, we write
\begin{equation}
D^M_{end} \: = \: 
\begin{pmatrix}
{\cal A} & {\cal B} \\
{\cal C} & {\cal D}
\end{pmatrix}.
\end{equation}

The operators ${\cal B}$ and
${\cal C}$ are bounded, as can be seen by the method of proof of
\cite[Proposition 2]{Lott2 (1999)},
replacing the operator $\widehat{d} \: + \: \widehat{d}^*$ of
\cite[Pf. of Proposition 2]{Lott2 (1999)} by $D^{Z_I}$.
As in \cite[Proposition 3]{Lott2 (1999)}, the
operator ${\cal D}$ has vanishing essential spectrum.
Put ${\cal L} \: = \: \begin{pmatrix}
{\cal A} & 0 \\
0 & {\cal D}
\end{pmatrix}$. To prove the theorem, it suffices to show that
$D^M_{end}$ and ${\cal L}$ have the same essential spectrum.
For this, it suffices to show that
$\left( D^M_{end} \: + \: k \: i \right)^{-1} \: - \:
\left( {\cal L} \: + \: k \: i \right)^{-1}$ is compact
for some $k \: > \: 0$
\cite[Vol. IV, Chapter XIII.4, Corollary 1]{Reed-Simon (1978)}.
 
We use the general identity that
\begin{equation} \label{inverse}
\begin{pmatrix}
\alpha & \beta \\
\gamma & \delta
\end{pmatrix}^{-1} \: = \:
\begin{pmatrix}
\alpha^{-1} \: + \: \alpha^{-1} \: \beta \: 
\left( \delta \: - \: \gamma \: \alpha^{-1} \: \beta \right)^{-1} \:
\gamma \: \alpha^{-1} & 
- \: \alpha^{-1} \: \beta \: 
\left( \delta \: - \: \gamma \: \alpha^{-1} \: \beta \right)^{-1} \\
- \: \left( \delta \: - \: \gamma \: \alpha^{-1} \: \beta \right)^{-1} \:
\gamma \: \alpha^{-1} &
\left( \delta \: - \: \gamma \: \alpha^{-1} \: \beta \right)^{-1}
\end{pmatrix}
\end{equation}
provided that $\alpha$ and 
$\delta \: - \: \gamma \: \alpha^{-1} \: \beta$ are invertible.  
Put 
\begin{equation}
\begin{pmatrix}
\alpha & \beta \\
\gamma & \delta
\end{pmatrix} \: = \: D^M_{end} \: + \: k \: i \: = \:
\begin{pmatrix}
{\cal A} \: + \: k \: i & {\cal B} \\
{\cal C} & {\cal D} \: + \: k \: i
\end{pmatrix}.
\end{equation}
If $k$ is positive then $\alpha$ and $\delta$ are invertible, with
$\delta^{-1}$ being compact.
If $k$ is large enough then
$\parallel
\delta^{- \: 1/2} \: \gamma \: \alpha^{-1} \: \beta \: 
\delta^{- \: 1/2} \parallel
\: < \: 1$.
Writing
\begin{equation} \label{Neumann}
\delta \: - \: \gamma \: \alpha^{-1} \: \beta \: = \:
\delta^{1/2} \: \left( I \: - \:
\delta^{- \: 1/2} \: \gamma \: \alpha^{-1} \: \beta
\: \delta^{- \: 1/2} \right)
\: \delta^{1/2},
\end{equation}
we now see that $\delta \: - \: \gamma \: \alpha^{-1} \: \beta$ is
invertible. It also follows from (\ref{Neumann}) that
$\left( \delta \: - \: \gamma \: \alpha^{-1} \: \beta \right)^{-1}$ is compact.
Using (\ref{inverse}), the theorem follows.

\end{document}